\documentclass[
    ,final            
  ]
  {aipproc}

\layoutstyle{8x11single}

  \usepackage{latexsym}
  \usepackage{amsmath}   
  \usepackage{amsfonts}
  \usepackage{amssymb}
  \usepackage{graphicx}
  \usepackage{amsthm}

  \newcommand{\field}[1]{\mathbb{#1}}
  
  \newcommand{\R}{\field{R}}

  \newcommand{\HQ}{\field{H}}
  \newcommand{\vect}[1]{\ensuremath{\mbox{\textbf{\textit{#1}}}}}
  \newcommand{\svect}[1]{\ensuremath{\mbox{\textbf{\textit{\small #1}}}}}
  \newcommand{\bomega}{\boldsymbol{\omega}}
  
  \newcommand{\qi}{\ensuremath{\mbox{\boldmath $i$}}}
  
  \newcommand{\qj}{\ensuremath{\mbox{\boldmath $j$}}}
  
  \newcommand{\qk}{\ensuremath{\mbox{\boldmath $k$}}}
  
  \newcommand{\be}{\begin{equation}}
  \newcommand{\ee}{\end{equation}} 
  
  \newcommand{\bfr}{\begin{frame}}
  \newcommand{\efr}{\end{frame}}

  \newtheorem{theorem}{Theorem}
  \newtheorem{lemma}{Lemma}

\begin{document}

\title{OPS-QFTs: A new type of quaternion Fourier transforms based on the orthogonal planes split with one or two general pure quaternions}

\classification{AMS Subj. Class. 15A66, 42A38}
\keywords      {Clifford geometric algebra, quaternion geometry, quaternion Fourier transform, inverse Fourier transform, orthogonal planes split}

\author{Eckhard Hitzer}{
  address={Department of Applied Physics, University of Fukui, 910-8507 Japan}
}

\begin{abstract}
We explain the orthogonal planes split (OPS) of quaternions based on the arbitrary choice of one or two linearly independent pure unit quaternions $f,g$.
Next we systematically generalize the quaternionic Fourier transform (QFT) applied to quaternion fields to conform with the OPS determined by $f,g$, or by only one pure unit quaternion $f$, comment on their geometric meaning, and establish inverse transformations.

\end{abstract}

\maketitle



References \cite{TAE:QFT,EH:QFTgen,EH:DirUP_QFT,SJS:FTcolQuat} give background on quaternion Fourier transformations. For details and proofs of the orthogonal planes split we refer the reader to \cite{HS:OPSofQ}.

Gauss, Rodrigues and Hamilton's four-dimensional (4D) quaternion algebra $\HQ$ is defined over $\R$ with three imaginary units:

\be
 \qi \qj = -\qj \qi = \qk, \,\,
 \qj \qk = -\qk \qj = \qi, \,\,
 \qk \qi = -\qi \qk = \qj, \,\, 
 \qi^2=\qj^2=\qk^2=\qi \qj \qk = -1.
\label{eq:quat}
\end{equation}
Every quaternion can be written explicitly as
\be
  q=q_r + q_i \qi + q_j \qj + q_k \qk \in \HQ, \quad 
  q_r,q_i, q_j, q_k \in \R,
  \label{eq:aquat}
\end{equation}
and has a \textit{quaternion conjugate} (equivalent to reversion in $Cl_{3,0}^+$)
\be
  \tilde{q} = q_r - q_i \qi - q_j \qj - q_k \qk,
  \quad 
  \widetilde{pq} = \tilde{q}\tilde{p},
  \label{eq:quatconj}
\end{equation} 
which leaves the scalar part $q_r$ unchanged.
This leads to the \textit{norm} of $q\in\HQ$
\be
  | q | = \sqrt{q\tilde{q}} = \sqrt{q_r^2+q_i^2+q_j^2+q_k^2},
  \qquad
  | p q | = | p || q |.
\end{equation}
The part $\vect{q} = q - q_r = (q-\tilde{q})/2 = q_i \qi + q_j \qj + q_k \qk$ is called a 
pure quaternion, and %
it squares to the negative number $-(q_i^2+q_j^2+q_k^2)$.
Every unit quaternion (i.e. $|q|=1$) can be written as 
\begin{gather}
  q 
  = q_r + q_i \qi + q_j \qj + q_k \qk
  = q_r + {\sqrt{q_i^2+q_r^2+q_k^2}}\,\hat{\vect{q}}
  = \cos \alpha + \sin \alpha \,\hat{\vect{q}} 
  = e^{\alpha \,\hat{\svect{q}}},
  \nonumber \\
  \cos \alpha = q_r, \quad 
  \sin \alpha = \sqrt{q_i^2+q_r^2+q_k^2}, \quad
  \hat{\vect{q}} = \frac{q_i \qi + q_j \qj + q_k\qk}{\sqrt{q_i^2+q_r^2+q_k^2}},
  \quad
  \hat{\vect{q}}^2 = -1.
\end{gather}
The inverse of a non-zero quaternion is
\be 
  q^{-1} = \frac{\tilde{q}}{|q|^2} = \frac{\tilde{q}}{q\tilde{q}}.
\ee 
The symmetric scalar part of a quaternion is defined as
\be
  Sc(q) = q_r = \frac{1}{2}(q+\tilde{q}),
  \qquad 
  Sc(pq) = Sc(qp) = p_rq_r - p_iq_i - p_jq_j - p_kq_k,
\ee
with linearity
\be 
  Sc(\alpha p+ \beta q) 
  = \alpha Sc(p) + \beta Sc(q) 
  = \alpha p_r + \beta q_r,
  \quad 
  \forall p,q \in \HQ, \,\,\, \alpha, \beta \in \R.
  \label{eq:Sclin}
\ee 
The scalar part and the quaternion conjugate allow the definition of 
the $\R^4$ \textit{inner product} of two quaternions $p,q$ as
\be 
  Sc(p\widetilde{q}) 
  = p_rq_r + p_iq_i + p_jq_j + p_kq_k \in \R . 
  \label{eq:4Dinnp}
\ee


We consider an arbitrary pair of linearly independent nonorthogonal pure quaternions $f,g$, $f^2=g^2=-1, f\neq \pm g$. The orthogonal 2D planes split (OPS) is then defined with respect to the linearly independent pure unit quaternions $f, g$ as 
\be 
q_{\pm} = \frac{1}{2}(q \pm f q g).
\ee 
We thus observe, that $f q g = q_+ - q_-$, i.e. under the map $f( )g$ the $q_+$ part is invariant, but the $q_-$ part changes sign.

Both parts are two-dimensional, and span two completely orthogonal planes. The $q_+$ plane is spanned by the orthogonal quaternions $\{f-g, 1+fg\}$
whereas the  $q_-$ plane is e.g. spanned by $\{f+g, 1-fg\}$. 

\begin{lemma}[Orthogonality of two OPS planes]\label{lm:OPSortho}
  Given two quaternions $q,p$ and applying the OPS with respect to two linearly independent pure unit quaternions $f,g$ we get zero for the
  scalar part of the mixed products
  \be
     Sc(p_+\widetilde{q}_-) = 0, \qquad Sc(p_-\widetilde{q}_+) = 0 .
  \ee
\end{lemma}
The set $\{f-g, 1+fg, f+g, 1-fg\}$ forms an orthogonal basis of $\HQ$ interpreted as $\R^4$. We can therefore use the following representation for every $q \in \HQ$ by means of four real coefficients $q_1, q_2, q_3, q_4 \in \R$
\begin{gather} 
  q = q_1 (1+fg) + q_2 (f-g) + q_3 (1-fg) + q_4 (f+g),
  \\
  q_1 = Sc(q(1+fg)^{-1}), \quad 
  q_2 = Sc(q(f-g)^{-1}), \quad
  q_3 = Sc(q(1-fg)^{-1}), \quad
  q_4 = Sc(q(f+g)^{-1}).
  \nonumber 
\end{gather}
As an example we have for $f=\qi, g=\qj$ we obtain
\be 
  q_1 = \frac{1}{2}(q_r+q_k), \quad 
  q_2 = \frac{1}{2}(q_i-q_j), \quad
  q_3 = \frac{1}{2}(q_r-q_k), \quad
  q_4 = \frac{1}{2}(q_i+q_j).
\ee 

\begin{theorem}[Determination of $f,g$ from given analysis planes]\label{th:detfg}
Assume that two desired analysis planes are given by a set of four orthogonal quaternions $\{a, b\}$ and $\{c, d\}$. Without restriction of generality $a^2=c^2=-1$. 

For the $\{a, b\}$ plane to become the $q_-$ plane and the  $\{c, d\}$ plane the $q_+$ plane we need to set
\be 
  f= ba, \qquad g = Sc(f\tilde{a}) a - Sc(f\tilde{c}) c.
\ee 
For the opposite assignment of the $\{a, b\}$ plane to become the $q_+$ plane and the  $\{c, d\}$ plane the $q_-$ plane, we only need to change the sign of $g$, i.e. we can set
\be 
  f= ba, \qquad g = Sc(f\tilde{c}) c - Sc(f\tilde{a}) a.
\ee 
\end{theorem}

The map $f( )g$ rotates the $q_-$ plane by $180^{\circ}$ around the $q_+$ axis plane. This interpretation of the map $f( )g$ is in perfect agreement with Coxeter's notion of \textit{half-turn} in \cite{HSMC:QuatRef}. 

The following identities hold
\be 
  e^{\alpha f} q_{\pm} e^{\beta g} 
  = q_{\pm} e^{(\beta\mp\alpha) g}
  = e^{(\alpha\mp\beta) f}q_{\pm}.
  \label{eq:eqpmeres}
\ee


The \textit{general} double sided \textit{orthogonal planes} (i.e. 2D subspaces) \textit{split quaternion Fourier transform} (OPS-QFT) is defined as
\begin{equation}
  \label{eq:QFTfg}
  \mathcal{F}^{f,g}\{ h \}(\bomega)
  = \int_{\R^2} e^{-f x_1\omega_1} h(\vect{x}) \,e^{-g x_2\omega_2} d^2\vect{x},
\end{equation}
where $h\in L^1(\R^2,\mathbb{H})$, $d^2\vect{x} = dx_1dx_2$ and $\vect{x}, \bomega \in \R^2 $. The OPS-QFT \eqref{eq:QFTfg} is invertible
\begin{equation}
  \label{eq:QFTfginv}
  h(\vect{x})
  = \frac{1}{(2\pi)^2}\int_{\R^2} e^{f x_1\omega_1}  \mathcal{F}^{f,g}\{ h \}(\bomega) \,e^{g x_2\omega_2} d^2\bomega,
\end{equation}

Linearity of the integral allows us to use the OPS split $h=h_- + h_+$, 
$h_{\pm}=\frac{1}{2}(h\pm fhg)$
\be 
  \mathcal{F}^{f,g}\{ h \}(\bomega)
  = \mathcal{F}^{f,g}\{ h_- \}(\bomega) + \mathcal{F}^{f,g}\{ h_+ \}(\bomega)
  = \mathcal{F}^{f,g}_-\{ h \}(\bomega) + \mathcal{F}_+^{f,g}\{ h \}(\bomega),
\ee 
since by its construction the operators of the Fourier transformation $\mathcal{F}^{f,g}$, and of the OPS with respect to $f,g$ commute. 
From \eqref{eq:eqpmeres} follows 
\begin{theorem}[OPS-QFT of $h_{\pm}$]
\label{th:fpmtrafo}
The QFT of the $h_{\pm}$ OPS split parts, with respect to two linearly independent unit quaternions $f,g$, of a quaternion module function 
$h \in L^2(\R^2,\HQ)$ have the complex forms
\be
  \mathcal{F}^{f,g}_{\pm}\{h\} 
  \stackrel{}{=} \int_{\R^2}
    h_{\pm}e^{-g (x_2\omega_2 \mp x_1\omega_1)}d^2x
  \stackrel{}{=} \int_{\R^2}
    e^{-f (x_1\omega_1 \mp x_2\omega_2)}h_{\pm}d^2x \,\, .
\end{equation}
\end{theorem}

The geometric interpretation of the integrand
  $e^{-f x_1\omega_1} h(\vect{x}) \,e^{-g x_2\omega_2}$
of the QFT${}^{f,g}$ in \eqref{eq:QFTfg} is: The integrand means to locally rotate by the phase angle $-(x_1\omega_1+x_2\omega_2)$ in the $q_-$ plane, and by phase angle  $-(x_1\omega_1-x_2\omega_2) = x_2\omega_2-x_1\omega_1$ in the $q_+$ plane.


The \textit{phase angle} OPS-QFT with a straight forward two phase angle interpretation is
\begin{equation}
  \label{eq:QFT2angle}
  \mathcal{F}_D^{f,g}\{ h \}(\bomega)
  = \int_{\R^2} e^{-f \frac{1}{2}(x_1\omega_1+x_2\omega_2)} h(\vect{x}) 
    \,e^{-g \frac{1}{2}(x_1\omega_1-x_2\omega_2)} d^2\vect{x}.
\end{equation}
where again $h\in L^1(\R^2,\mathbb{H})$, $d^2\vect{x} = dx_1dx_2$ and $\vect{x}, \bomega \in \R^2 $. Its inverse is given by
\begin{equation}
  \label{eq:QFT2angleinv}
  h(\vect{x}) 
  = \frac{1}{(2\pi)^2}\int_{\R^2} e^{f \frac{1}{2}(x_1\omega_1+x_2\omega_2)} 
    \mathcal{F}_D^{f,g}\{ h \}(\bomega)
    \,e^{g \frac{1}{2}(x_1\omega_1-x_2\omega_2)} d^2\bomega.
\end{equation}

The geometric interpretation of the integrand of \eqref{eq:QFT2angle} is a local phase rotation by angle 
$-(x_1\omega_1+x_2\omega_2)/2- (x_1\omega_1-x_2\omega_2)/2 = - x_1\omega_1$ 
in the $q_-$ plane, and a second local phase rotation by angle 
$-(x_1\omega_1+x_2\omega_2)/2+ (x_1\omega_1-x_2\omega_2)/2 = - x_2\omega_2$
in the $q_+$ plane. 

If we apply the OPS${}^{f,g}$ split to \eqref{eq:QFT2angle} we obtain the following two parts
\be 
  \mathcal{F}^{f,g}_{D+}\{h\} 
  \stackrel{}{=} \int_{\R^2}
    h_{+}e^{+g x_2\omega_2 }d^2x
  \stackrel{}{=} \int_{\R^2}
    e^{-f x_2\omega_2}h_{+}d^2x 
  \, , \quad
  \mathcal{F}^{f,g}_{D-}\{h\} 
  \stackrel{}{=} \int_{\R^2}
    h_{-}e^{-g x_1\omega_1 }d^2x
  \stackrel{}{=} \int_{\R^2}
    e^{-f x_1\omega_1}h_{-}d^2x .
\end{equation}


The OPS with respect to, e.g., $f=g=\qi$ gives
\be 
  q_{\pm} = \frac{1}{2}(q\pm \qi q \qi), \quad 
  q_+ = q_j \qj + q_k \qk = (q_j + q_k \qi)\qj, \quad 
  q_- = q_r + q_i \qi ,
  \label{eq:opsiqi}
\ee  
where the $q_+$ plane is two-dimensional and manifestly orthogonal to the 2D $q_-$ plane. The above corresponds to the simplex/perplex split of \cite{TAE:QFT}. $e^{\alpha f} q\, e^{\beta f}$ means a rotation by angle $\alpha+\beta$ in the $q_-$ plane followed by a rotation by angle $\alpha-\beta$ in the orthogonal $q_+$ plane.

A variant of the OPS-QFT with $g=f$ is therefore
\begin{equation}
  \label{eq:QFTff}
  \mathcal{F}^{f,f}\{ h \}(\bomega)
  = \int_{\R^2} e^{-f x_1\omega_1} h(\vect{x}) \,e^{-f x_2\omega_2} d^2\vect{x},
\end{equation}
where $h\in L^1(\R^2,\mathbb{H})$, $d^2\vect{x} = dx_1dx_2$ and $\vect{x}, \bomega \in \R^2 $. The immediate geometric interpretation is that the integrand $e^{-f x_1\omega_1} h(\vect{x}) \,e^{-f x_2\omega_2}$ leads to a local phase roation by angle $-(x_1\omega_1+x_2\omega_2)$ in the $q_-$ plane combined with a second local phase rotation by angle $x_2\omega_2-x_1\omega_1$ in the $q_+$ plane. The inverse transform is given by 
\begin{equation}
  \label{eq:QFTffinv}
  h(\vect{x})
  = \frac{1}{(2\pi)^2} \int_{\R^2} e^{f x_1\omega_1} \mathcal{F}^{f,f}\{ h \}(\bomega)  \,e^{f x_2\omega_2} d^2\bomega,
\end{equation}

The \textit{phase angle} OPS-QFT, using $g=f$, with a straight forward two phase angle interpretation is
\begin{equation}
  \label{eq:QFT2angleff}
  \mathcal{F}_D^{f,f}\{ h \}(\bomega)
  = \int_{\R^2} e^{-f \frac{1}{2}(x_1\omega_1+x_2\omega_2)} h(\vect{x}) 
    \,e^{-f \frac{1}{2}(x_1\omega_1-x_2\omega_2)} d^2\vect{x}.
\end{equation}
where again $h\in L^1(\R^2,\mathbb{H})$, $d^2\vect{x} = dx_1dx_2$ and $\vect{x}, \bomega \in \R^2 $. Its inverse is given by
\begin{equation}
  \label{eq:QFT2angleffinv}
  h(\vect{x}) 
  = \frac{1}{(2\pi)^2}\int_{\R^2} e^{f \frac{1}{2}(x_1\omega_1+x_2\omega_2)} 
    \mathcal{F}_D^{f,f}\{ h \}(\bomega)
    \,e^{f \frac{1}{2}(x_1\omega_1-x_2\omega_2)} d^2\bomega.
\end{equation}

The geometric interpretation of the integrand of \eqref{eq:QFT2angleff} is a local phase rotation by angle 
$-(x_1\omega_1+x_2\omega_2)/2- (x_1\omega_1-x_2\omega_2)/2 = - x_1\omega_1$ 
in the $q_-$ plane, and a second local phase rotation by angle 
$-(x_1\omega_1+x_2\omega_2)/2+ (x_1\omega_1-x_2\omega_2)/2 = - x_2\omega_2$
in the $q_+$ plane. 

If we apply the OPS${}^{f,f}$ split to \eqref{eq:QFT2angleff} we obtain the following two parts
\be 
  \mathcal{F}^{f,f}_{D+}\{h\} 
  \stackrel{}{=} \int_{\R^2}
    h_{+}e^{+f x_2\omega_2 }d^2x
  \stackrel{}{=} \int_{\R^2}
    e^{-f x_2\omega_2}h_{+}d^2x 
  \, , \quad
  \mathcal{F}^{f,f}_{D-}\{h\} 
  \stackrel{}{=} \int_{\R^2}
    h_{-}e^{-f x_1\omega_1 }d^2x
  \stackrel{}{=} \int_{\R^2}
    e^{-f x_1\omega_1}h_{-}d^2x .
\end{equation}

The general OPS-QFT involving \textit{quaternion conjugation} is
defined as
\begin{equation}
  \label{eq:QFTgf}
  \mathcal{F}_c^{g,f}\{ h \}(\bomega)
  = \int_{\R^2} e^{-g x_1\omega_1} \widetilde{h(\vect{x})} \,e^{-f x_2\omega_2} d^2\vect{x},
\end{equation}
where $h\in L^1(\R^2,\mathbb{H})$, $d^2\vect{x} = dx_1dx_2$ and $\vect{x}, \bomega \in \R^2 $. The inverse is given by
\begin{equation}
  \label{eq:QFTgfinv}
  h(\vect{x})
  = \frac{1}{(2\pi)^2}\int_{\R^2} e^{-f x_2\omega_2} \widetilde{\mathcal{F}_c^{g,f}\{ h \}(\bomega)} \,e^{-g x_1\omega_1} d^2\bomega,
\end{equation}

This approach results in the following OPS theorem.
\begin{theorem}[OPS-QFT $\mathcal{F}_c^{g,f}$ of $h_{\pm}$]
\label{th:gfpmtrafo}
The OPS-QFT $\mathcal{F}_c^{g,f}$ \eqref{eq:QFTgf} of the $h_{\pm} = \frac{1}{2}(h\pm fhg)$ OPS split parts, with respect to two linearly independent unit quaternions $f,g$, of a quaternion module function 
$h \in L^2(\R^2,\HQ)$ have the complex forms
\be
  \mathcal{F}^{g,f}_{c,\pm}\{h\} 
  \stackrel{}{=} \int_{\R^2}
    \widetilde{h_{\pm}}e^{-f (x_2\omega_2 \mp x_1\omega_1)}d^2x
  \stackrel{}{=} \int_{\R^2}
    e^{-g (x_1\omega_1 \mp x_2\omega_2)}\widetilde{h_{\pm}}d^2x \,\, .
\end{equation}
\end{theorem}

The variant of the general OPS-QFT involving \textit{quaternion conjugation}, using $g=f$, is
defined as
\begin{equation}
  \label{eq:QFTffc}
  \mathcal{F}_c^{f,f}\{ h \}(\bomega)
  = \int_{\R^2} e^{-f x_1\omega_1} \widetilde{h(\vect{x})} \,e^{-f x_2\omega_2} d^2\vect{x},
\end{equation}
where $h\in L^1(\R^2,\mathbb{H})$, $d^2\vect{x} = dx_1dx_2$ and $\vect{x}, \bomega \in \R^2 $. The inverse is given by
\begin{equation}
  \label{eq:QFTffcinv}
  h(\vect{x})
  = \frac{1}{(2\pi)^2}\int_{\R^2} e^{-f x_2\omega_2} \widetilde{\mathcal{F}_c^{f,f}\{ h \}(\bomega)} \,e^{-f x_1\omega_1} d^2\bomega ..
\end{equation}

This results in the following OPS theorem.
\begin{theorem}[OPS-QFT $\mathcal{F}_c^{f,f}$ of $h_{\pm}$]
\label{th:ffpmtrafo}
The OPS-QFT $\mathcal{F}_c^{f,f}$ \eqref{eq:QFTffc} of the $h_{\pm} = \frac{1}{2}(h\pm fhf)$ OPS split parts, with respect to the unit quaternion $f$, of a quaternion module function 
$h \in L^2(\R^2,\HQ)$ have has complex forms
\be
  \mathcal{F}^{f,f}_{c,\pm}\{h\} 
  \stackrel{}{=} \int_{\R^2}
    \widetilde{h_{\pm}}e^{-f (x_2\omega_2 \mp x_1\omega_1)}d^2x
  \stackrel{}{=} \int_{\R^2}
    e^{-f (x_1\omega_1 \mp x_2\omega_2)}\widetilde{h_{\pm}}d^2x \,\, .
\end{equation}
\end{theorem}

\begin{theacknowledgments}
  Soli deo gloria. I do thank my dear family, S. Sangwine, W. Spr\"{o}ssig and K. G\"{u}rlebeck.
\end{theacknowledgments}

\bibliographystyle{aipproc}

\end{document}